\documentclass[10pt,twoside]{siamltex}
\usepackage{amsfonts,epsfig}
\usepackage{amssymb} 
\usepackage{amsmath}

\newtheorem{example}[theorem]{Example}

\begin{document}

\bibliographystyle{plain}
\title{On Determining Minimal Spectrally Arbitrary Patterns
\thanks{Received by the editors on ...
Accepted for publication on ....   Handling Editor: ...}}

\author{Michael S. Cavers\thanks{Department of Combinatorics and 
Optimization, University of Waterloo, Waterloo, ON Canada, 
N2L 3G1 (mscavers@math.uwaterloo.ca).}
\and
In-Jae Kim\thanks{Department of Mathematics, University of Wyoming, Laramie, WY 82071, USA
\newline (injaekim89@hotmail.com, bshader@uwyo.edu).}
\and
Bryan L. Shader\footnotemark[3]
\and
Kevin~N.~Vander~Meulen\thanks{Department of Mathematics, Redeemer University College, 
Ancaster, ON Canada, L9K 1J4 (kvanderm@redeemer.on.ca). Research supported in part by NSERC.}}


\pagestyle{myheadings}
\markboth{M.S.\ Cavers, I-J.\ Kim, B.L.\ Shader and K.N.\ Vander Meulen}{On Determining Spectrally Arbitrary Patterns}
\maketitle

\begin{abstract}
In this paper we present a new family of minimal spectrally arbitrary patterns which allow for arbitrary spectrum by using the Nilpotent-Jacobian method introduced in \cite{DJOD}.  The novel approach here is that we use the Intermediate Value Theorem to avoid finding an explicit nilpotent realization of the new minimal spectrally arbitrary patterns.
\end{abstract}

\begin{keywords}
sign pattern, spectrum, nilpotent matrix 
\end{keywords}
\begin{AMS}
15A18, 15A29
\end{AMS}

\section{Introduction} \label{intro-sec}
A matrix $\mathcal{S}$ with entries in $\{+,-,0\}$ is a \emph{sign pattern}.
Let $\mathcal{S}=[s_{ij}]$ and $\mathcal{U}=[u_{ij}]$ be $m$ by $n$ sign patterns.  
If $u_{ij}=s_{ij}$ whenever $s_{ij} \ne 0$, 
then $\mathcal{U}$ is a \emph{superpattern} of $\mathcal{S}$ and $\mathcal{S}$ is a \emph{subpattern} of $\mathcal{U}$. 
A subpattern of $\mathcal{S}$ which is not $\mathcal{S}$ itself is a \emph{proper subpattern} of $\mathcal{S}$.
Similarly, a superpattern of $\mathcal{S}$ which is not $\mathcal{S}$ itself is a \emph{proper superpattern} of $\mathcal{S}$.

For a given real number $a$, the sign of $a$ is denoted by 
$\mbox{sgn}(a)$, and is $+$, $0$, or $-$ depending on whether $a$ is positive, $0$ or negative.   
The \emph{sign pattern class} of an $m$ by $n$ sign pattern $\mathcal{S}=[s_{ij}]$ is defined by

\begin{center} 
$Q(\mathcal{S})=\{ A=[a_{ij}] \in \mathbb{R}^{m \times n}$ : $\mbox{sgn}(a_{ij})=s_{ij}$ for all $i, j\}$.  
\end{center}

\noindent
A \emph{nilpotent realization} of a sign pattern $\mathcal{S}$ is a real matrix
$A\in Q(\mathcal{S})$ whose only eigenvalue is zero. We say a sign pattern $\mathcal{S}$ {\emph{requires}} a property $P$ if
each matrix $A\in Q(\mathcal{S})$  has property $P$.

For a complex number $\lambda$, we take $\overline{\lambda}$ to denote the complex conjugate of $\lambda$.
Let $\sigma$ be a multi-list of complex numbers.
Then $\sigma$ is \emph{self-conjugate} if and only if for each $\lambda \in \sigma$, $\overline{\lambda}$ occurs with
the same multiplicity as $\lambda$'s in $\sigma$.
Note that if $A$ is an $n$ by $n$ real matrix, then the spectrum of $A$ is self-conjugate. 

An $n$ by $n$ sign pattern $\mathcal{S}$ is a \emph{spectrally arbitrary pattern} (\emph{SAP}) if each self-conjugate multi-list
of $n$ complex numbers is the spectrum of a realization of $\mathcal{S}$, that is,
if each monic real polynomial of degree $n$ is the characteristic polynomial of a matrix in $Q(\mathcal{S})$.
If $\mathcal{S}$ is a SAP and no proper subpattern of  $A$ is spectrally arbitrary, then $\mathcal{S}$ is a \emph{minimal spectrally arbitrary pattern} (\emph{MSAP}). 

The question of the existence of a SAP arose in \cite{DJOD}, 
where a general method (based on the Implicit Function Theorem) was given to  
prove that a sign pattern and all of its superpatterns are SAP.
The first SAP of order $n$ for each $n \geq 2$ was provided in \cite{MOTD}.
Later, the method in \cite{DJOD}, which we will call the  \emph{Nilpotent-Jacobian method $($N-J method$)$}, 
was reformulated in \cite{BMOD}, and used to find MSAPs 
in \cite{BMOD},~\cite{CV} and \cite{MTD}.

The N-J method is quite powerful, yet its Achilles' heel is 
the need to determine (not necessarily explicitly) an appropriate
nilpotent realization 
in order to compute the Jacobian involved in the method.
It is not an easy task to find an appropriate nilpotent realization, even for the well-structured antipodal tridiagonal pattern 
in \cite{DJOD}, at which the Jacobian is nonzero (see \cite{EOD}).
In this paper, we show how to use the N-J method without explicitly constructing a nilpotent realization 
by using the Intermediate Value Theorem, and 
provide a new family of MSAPs.

\section{The N-J Method}\label{NJ}

Throughout, we take $p_{A}(x)=x^{n}- \alpha_{1}x^{n-1}+\alpha_{2}x^{n-2}-\cdots+(-1)^{n-1}\alpha_{n-1}x+(-1)^{n}\alpha_{n}$ 
to denote the characteristic polynomial of an $n$ by $n$ matrix $A$, and 
the Jacobian $\Delta=\det \left( \displaystyle\frac{\partial f_{i}}{\partial x_{j}} \right)$ is denoted by
$\displaystyle\frac{\partial (f_1, \ldots, f_n)}{\partial (x_1 , \ldots, x_n)}$ where
$f=(f_{1}, \ldots, f_{n})$ is a function of $x_{1}, \ldots, x_{n}$ such that $\displaystyle\frac{\partial f_{i}}{\partial x_{j}}$ 
exists for all $i, j \in \{1, \ldots, n\}$.
The matrix $\left( \displaystyle\frac{\partial f_{i}}{\partial x_{j}} \right)$ is called the \emph{Jacobian matrix}
of $f$.

The next theorem describes the N-J method for proving that
a sign pattern and all of its superpatterns are spectrally arbitrary.

\medskip

\begin{theorem}{\mbox{\rm{(\cite{BMOD}, Lemma 2.1)}}}
\label{theoremNJ}
Let $\mathcal{S}$ be an $n$ by $n$ sign pattern, and suppose that there exists a nilpotent realization $M=[m_{ij}]$ of $\mathcal{S}$ 
with at least $n$ nonzero entries, say, $m_{i_1 j_1}, \ldots, m_{i_n j_n}$.  
Let $X$ be the matrix obtained by replacing these entries in $M$ by variables $x_1 , \ldots, x_n$, and let 
$(-1)^{k}\alpha_{k}$ be the coefficients of $p_{X}(x)$ for $k=1, 2, \ldots, n$.
If the Jacobian $\displaystyle\frac{\partial (\alpha_1, \ldots, \alpha_n)}{\partial (x_1 , \ldots, x_n)}$ is nonzero at $(x_1 , \ldots, x_n )=(m_{i_1 j_1}, \ldots, m_{i_n j_n})$, then every superpattern of $\mathcal{S}$ is spectrally arbitrary.
\end{theorem}

\begin{example}{\mbox{\rm{(\cite{BMOD}, Example $2.2$)}}}
\label{ex:2}
Let $\mathcal{S} = \left[ \begin{array}{cc} + & - \\  + & -  \end{array} \right]$.
Then $A = \left[ \begin{array}{cc} 1 & -1 \\  1 & -1  \end{array} \right]$ is a nilpotent realization of $\mathcal{S}$.
Let $X= \left[ \begin{array}{cr} x_{1} & -1 \\  1 & x_{2}  \end{array} \right]$. 
Then,
$$
p_{X}(x)= x^{2} - \alpha_{1}x + \alpha_{2},
$$
where $\alpha_{1}=x_{1}+x_{2}$ and $\alpha_{2} = x_{1}x_{2}+1$.
Thus
$$
\Delta= \displaystyle\frac{\partial (\alpha_1,\alpha_{2})}{\partial (x_1 ,x_2)} = \det \left[ \begin{array}{cc} 1 & 1 \\  x_{2} & x_{1} \end{array} \right]
=x_{1} - x_{2}.
$$
At $(x_{1}, x_{2}) = (1, -1)$, $\Delta=2 \ne 0$.
By Theorem~$\ref{theoremNJ}$, $\mathcal{S}$ is spectrally arbitrary.

If some entries of $\mathcal{S}$ are replaced by $0$, then $\alpha_{2}$ of each realization of the resulting sign pattern has a fixed sign 
in $\{+, -,0 \}$ and hence,
the resulting sign pattern is not spectrally arbitrary.
Therefore, $\mathcal{S}$ is a MSAP.
\end{example}

\section{A new MSAP}\label{new}

Let $n$ and $r$ be positive integers with $2 \leq r \leq n$, and let 
$$
\mathcal{K}_{n,r}=
\left[ \begin{array}{cccccccc} 
+      & -      & 0      & 0      & \cdots & \cdots & \cdots & 0      \\
+      & 0      & -      & 0      &        &        &        & \vdots \\
+     & 0      & 0      & -      & \ddots &        &        & \vdots  \\ 
\vdots & \vdots &        & \ddots & \ddots & \ddots &        & \vdots  \\
       &        &        &        & \ddots & \ddots & 0      & 0        \\
\vdots & \vdots &        &        &        & \ddots & -      & 0        \\
+      & 0      & \cdots &        &        & \cdots & 0      & -        \\
0      & \cdots & 0      & +      & 0      & \cdots & 0      & -        
\end{array} \right]_{n \times n},
$$

\noindent
where the positive entry in the last row is in column $n-r+1$.
For $r > \frac{n}{2}$, the patterns $\mathcal{K}_{n,r}$ are  spectrally arbitrary  \cite[Theorem 4.4]{CV}.
The pattern $\mathcal{K}_{n,n}$ was shown to be a MSAP in~\cite{BMOD} and
$\mathcal{K}_{n,n-1}$ was shown to be a MSAP in~\cite{CV}.
In this paper, we show that $\mathcal{K}_{n,r}$ is a SAP (and in fact a MSAP) 
for all $r$ with $2\leq r< n$.

It is convenient to consider matrices $A \in Q(\mathcal{K}_{n,r})$ of the form 
\begin{equation}
\label{matrix}
A=
\left[ \begin{array}{cccccccc} 
a_{1}      & -1      & 0      & 0      & \cdots & \cdots & \cdots & 0      \\
a_{2}      & 0       & -1     & 0      &        &        &        & \vdots \\
a_{3}     & 0      & 0       & -1     & \ddots &        &        & \vdots  \\ 
\vdots     & \vdots &         & \ddots & \ddots & \ddots &        & \vdots  \\
           &        &         &        & \ddots & \ddots & 0      & 0        \\
\vdots     & \vdots &         &        &        & \ddots & -1     & 0        \\
a_{n-1}    & 0      & \cdots  & 0      &\cdots        & \cdots & 0      & -1       \\
0           & \cdots & 0      & b_{r}  & 0      & \cdots & 0      & -1        
\end{array} \right]_{n \times n},
\end{equation}

\noindent
where $a_{j} > 0$ for $j=1, \ldots, n-1$ and $b_{r} >0$. 

We first give a definition and a little result on the positive
 zeros of real polynomials (of finite degree).
For a real polynomial $f(t)$, we set
$$
Z_{f} =\{a > 0~ | ~f(a)=0\}.
$$
\noindent
If $Z_{f}$ is nonempty, the minimum of $Z_{f}$ is denoted by $\min(Z_{f})$.

\medskip

\begin{proposition}
\label{IVT}
Let $f(t)$ and $g(t)$ be real polynomials, and $h(t)=g(t)-tf(t)$.
Suppose that $f(0)$, $g(0) > 0$, $Z_{f}$ and $Z_{g}$ are nonempty, and $\min(Z_{g}) < \min(Z_{f})$.
Then $Z_{h}$ is nonempty and $\min(Z_{h}) < \min(Z_{g})$.
\end{proposition}

\begin{proof}
Let $t_{g} = \min(Z_{g})$.  
Since $t_g < \min(Z_{f})$ and $f(0)>0$, we have $f(t_g) >0$.  
Thus,
\begin{center}
$h(t_g)=g(t_g)-t_g f(t_g)= -t_g f(t_g) < 0$.
\end{center}

Since $h(0)=g(0)>0$, the Intermediate Value Theorem implies that there exists a real number $a$ such that $0 < a < t_g$ and $h(a)=0$.
Thus, $Z_h$ is nonempty and $\min (Z_{h}) < \min (Z_{g})$. 
\end{proof}

Using Proposition~\ref{IVT}, we show the existence of a nilpotent realization of $\mathcal{K}_{n,r}$.

\medskip

\begin{lemma}
\label{NP}
For $n \geq 3$, $\mathcal{K}_{n,r}$ has a nilpotent realization.
\end{lemma}

\begin{proof}
Let $A \in Q(\mathcal{K}_{n,r})$ be of the form~(\ref{matrix}).
For convenience, we set $a_{0}=1$.
From~\cite[p. 262]{BMOD}, we can deduce that the coefficients of the characteristic polynomial
of $A$ satisfy
\begin{equation}
\label{char}
\begin{tabular}{rcl}
$\alpha_{1}$ & $=$ & $a_{1}-1$, \\
$\alpha_{j}$ & $=$ & $a_{j}-a_{j-1}$ \mbox{for} $j=2, \ldots, r-1$, \\
$\alpha_{j}$ & $=$ & $a_{j}-a_{j-1}+b_{r}a_{j-r}$ \mbox{for} $j=r, \ldots, n-1$, \mbox{and}  \\
$\alpha_{n}$ & $=$ & $b_{r}a_{n-r}-a_{n-1}$.
\end{tabular}
\end{equation}

Let $a_{1}=\cdots=a_{r-1}=1$, $b_{r}=t$.
Then $\alpha_{1}= \alpha_{2} = \cdots = \alpha_{r-1} =0$.
Note that if $n=r$,  then  setting $t=1$ gives a solution to (\ref{char}), and 
hence $\mathcal{K}_{n,n}$ has a nilpotent realization.
Suppose $n>r$.
In order to show that there exists a nilpotent realization of $\mathcal{K}_{n,r}$, we need to show the existence of positive numbers 
$a_{r},a_{r+1}, \ldots, a_{n-1}, t$ satisfying the following equations (obtained by setting $\alpha_{j}$'s to be zero for all $j=r, \ldots, n$):
\begin{equation}
\label{polys}
\begin{tabular}{rcl}
$a_{r}(t)$        & = & $a_{r-1}(t)-ta_{0}(t) = 1-t$ \\
$a_{r+1}(t)$      & = & $a_{r}(t)-ta_{1}(t) = 1-2t$ \\
$a_{r+2}(t)$      & = & $a_{r+1}(t)-ta_{2}(t) = 1-3t$ \\
                  & $\vdots$ &  \\
$a_{2r-1}(t)$     & = & $a_{2r-2}(t) -ta_{r-1}(t)= 1-rt $\\
$a_{2r}(t)$       & = & $a_{2r-1}(t) -t a_{r}(t)$ \\
$a_{2r+1}(t)$     & = & $a_{2r}(t) -ta_{r+1}(t)$ \\
                  &$\vdots$&  \\
$a_{n-1}(t)$      & = & $a_{n-2}(t) -ta_{n-1-r}(t)$ \\                  
$0$               & = & $a_{n-1}(t)-t a_{n-r}(t)$.
\end{tabular}
\end{equation}

\noindent For $j=r,r+1, \ldots,n-1$,
the polynomials in~(\ref{polys}) satisfy 
\begin{equation}
\label{poly}
a_{j}(t)=a_{j-1}(t)-ta_{j-r}(t)
\end{equation}

Let $h(t)=a_{n-1}(t)-t a_{n-r}(t)$.
It can be easily checked that $a_{j}(0)=1$ for all $j=r, r+1, \ldots, n-1$, and $h(0)=1$.
Since $a_{r+i}(t)=1-(i+1)t$ for $i=0,1, \ldots, r-1$, the zero of $a_{r+i}$ is $\displaystyle\frac{1}{i+1}$ for $i=0,1, \ldots, r-1$.
Hence,
\begin{equation}
\label{minZ}
\min(Z_{a_{2r-1}}) < \min(Z_{a_{2r-2}})< \cdots < \min(Z_{a_{r+1}}) < \min(Z_{a_{r}}).
\end{equation}

Since, by~(\ref{minZ}) $\min(Z_{a_{2r-1}}) < \min(Z_{a_{r}})$, and $a_{2r}(t)=a_{2r-1}(t)-ta_{r}(t)$, Proposition~\ref{IVT} implies that
$\min(Z_{a_{2r}}) < \min(Z_{a_{2r-1}})$.
Likewise, by repeatedly using Proposition~\ref{IVT}, we have
$$
\min(Z_{h}) < \min(Z_{a_{n-1}})< \cdots < \min(Z_{a_{r+1}}) < \min(Z_{a_{r}}).
$$

Thus, if $\min(Z_{h})$ is denoted by $t_{h}$, then $a_{j}(t_{h}) > 0$ for all $j=r, r+1, \ldots, n-1$.
Therefore, there exists a nilpotent realization of $\mathcal{K}_{n,r}$ for $n \geq 3$ 
when $a_{1}= \cdots = a_{r-1}=1$, $b_{r}=t_{h}$, and $a_{j}=a_{j}(t_{h})$ for $j=r, r+1, \ldots, n-1$.
\end{proof}

Throughout the remainder of this section $t_{h} = \min(Z_{h})$, and $a^{0}_{j}$ denotes the positive numbers $a_{j}(t_{h})$ 
for $j=r, r+1, \ldots, n-1$, and $a^{0}_{j}=1$ for $j=0,1, \ldots, r-1$, where $a_{j}(t)$'s are 
polynomials in~(\ref{polys}), and $h(t)=a_{n-1}(t)-t a_{n-r}(t)$.

While we made use of the characteristic polynomial associated with the 
patterns $\mathcal{V}^*_n(I)$ in~\cite{BMOD}, unlike for 
$\mathcal{V}^*_n(I)$,  for $\mathcal{K}_{n,r}$ we do not always end up with a 
Jacobian matrix whose pattern requires a signed determinant. 
 Thus we now develop some propositions to show 
that the Jacobian $\displaystyle\frac{\partial (\alpha_1, \ldots, \alpha_n)}{\partial (a_1 , \ldots, a_{n-1},b_{r})}$ is nonzero at the nilpotent realization of $\mathcal{K}_{n,r}$, $(1,a^{0}_{2}, \ldots, a^{0}_{n-1},t_h)$. 

First, we consider the matrix 
$$A_{k,r}=
\left[ \begin{array}{rrrcccrr} 
-1      & 1      & 0      & 0      & \cdots & \cdots & \cdots & 0      \\
0       & -1     & 1      & 0      &        &        &        & \vdots \\
\vdots  & 0      & -1     & 1      & \ddots &        &        & \vdots  \\ 
0       & \vdots &\ddots & \ddots & \ddots & \ddots &        & \vdots  \\
t_{h}   & 0      &        &        & \ddots & \ddots & 0      & 0        \\
0       & t_{h } & \ddots &        &        & \ddots & 1      & 0        \\
\vdots  & \ddots & \ddots &\ddots        &        & \ddots & -1     & 1       \\
0       & \cdots & 0      & t_{h}  & 0      & \cdots & 0      & -1        
\end{array} \right]_{k \times k},
$$
\vskip .1in
\noindent
where $t_{h}$ in the last row is in column $k-r+1$, and $k \geq 1$, $r \geq 2$.
If $r > k$, $A_{k,r}$ does not have any $t_{h}$ as an entry.
In addition, $\det(A_{0,r})$ is defined to be $1$.
Now, we find an explicit form of the determinant of $A_{k,r}$.

\medskip

\begin{proposition}
\label{detA}
If $2\leq r<n$ and $0 \leq k < n$ then
$$\det(A_{k,r})=
\begin{cases}
(-1)^{k}a^{0}_{k}& \quad when \quad 2 \leq r \leq k\\
(-1)^{k}& \quad when \quad  r > k
\end{cases}$$
\end{proposition}

\begin{proof}
The proof is by induction on $k$.
The cases for $k=0,1$ are clear.
Suppose $k=2$.
If $r=k$, then 
$$
\det(A_{2,2})=\det \left[ \begin{array}{cr} -1 & 1 \\  t_{h} & -1 \end{array} \right] =1-t_{h}.
$$
Since $r<n$, by~(\ref{polys}) we have $1-t_{h}=a^{0}_{r}$. 
Since $r=k=2$, $1-t_{h}=(-1)^{k}a^{0}_{k}$.
If $r > 2=k$, then
$$
\det(A_{2,r})=\det \left[ \begin{array}{rr} -1 & 1 \\ 0 & -1 \end{array} \right] = (-1)^{2}.
$$

Assume $k >2$.
For $2 < k \leq n-1$ and $r > k$, $A_{k,r}$ is an upper triangular matrix each of whose diagonal entries is $-1$.
Hence, $\det(A_{k,r})=(-1)^{k}$.

Suppose $r \leq k$.
By cofactor expansion along the first column, 
$$
\det(A_{k,r})=       (-1)\det(A_{k-1,r}) + (-1)^{r+1}t_{h}\det(A_{k-r,r}).
$$ 
If $k \geq 2r$, then  $k-r \geq r$ and $k-1 \geq r$.
By induction and~(\ref{poly}),
\begin{center}
\begin{tabular}{rcl}
$\det(A_{k,r})$& = & $(-1)(-1)^{k-1}a^{0}_{k-1}+ (-1)^{r+1}t_{h}(-1)^{k-r}a^{0}_{k-r}$\\
                                                     & = & $(-1)^{k}a^{0}_{k-1} + (-1)^{k+1}t_{h}a^{0}_{k-r}$\\
                                                     & = & $(-1)^{k}(a^{0}_{k-1}-t_{h}a^{0}_{k-r})$\\
                                                     & = &$(-1)^{k}a^{0}_{k}$ \quad since $k\geq 2r$.
\end{tabular}
\end{center}
\noindent
Next, suppose $r+1 \leq k < 2r$.
Since $k-1 \geq r$ and $k-r < r$, by induction, and~(\ref{poly}),
\begin{center}
\begin{tabular}{rcl}
$\det(A_{k,r})$& = & $(-1)(-1)^{k-1}a^{0}_{k-1}+ (-1)^{r+1}t_{h}(-1)^{k-r}$\\
                                                     & = & $(-1)^{k}a^{0}_{k-1} + (-1)^{k+1}t_{h}$\\
                                                     & = & $(-1)^{k}(a^{0}_{k-1}-t_{h})$\\ 
                                                     &=&$(-1)^{k}(a^{0}_{k-1}-t_{h}a^{0}_{k-r})$\\
                                                     &=&$(-1)^{k}a^{0}_{k}$.
\end{tabular}
\end{center}
\noindent
Lastly, suppose $r \leq k < r+1$, i.e. $k=r$.
By induction and~(\ref{polys}),
\begin{center}
\begin{tabular}{rcl}
$\det(A_{k,r})$& = & $(-1)(-1)^{k-1} + (-1)^{r+1}t_{h}(-1)^{k-r}$\\
                                                     & = & $(-1)^{k}(1-t_{h})$\\ 
                                                     &= &$(-1)^{k}a^{0}_{k}$.
\end{tabular}
\end{center}
\end{proof}

For $l \geq 1$, $r \geq 2$, and $c_{j} > 0$ ($j=1, \ldots, l$), let
$$
B_{l,r} = \left[ \begin{array}{rrrccccrc} 
1       & 0      & 0      & 0      & \cdots &\cdots     & \cdots & 0       & c_{1} \\
-1      & 1      & 0      & 0      &       &        &        & \vdots  & c_{2} \\
0       & -1     &  1     & 0      &        &        &        &         & \vdots \\ 
\vdots  &0 & -1     & 1      & \ddots &        &        &         & \vdots  \\
0       &  \vdots &\ddots    & \ddots & \ddots &        &        &         & \vdots \\
t_{h}   & 0      &        &        &           &        & \ddots & \vdots  & \vdots \\
0       & t_{h}  & \ddots &        &        &        & \ddots & 0       & \vdots \\
\vdots  & \ddots & \ddots &  \ddots     &        & \ddots & \ddots &  1      & c_{l-1} \\
0       & \cdots & 0      & t_{h}  & 0      & \cdots & 0      & -1      & c_{l}  
\end{array} \right]_{l \times l},
$$
where $t_{h}$ in the last row is in column $l-r$. 
If $r \geq l$, $B_{l,r}$ does not have any $t_{h}$ entry. 

\medskip

\begin{proposition}
\label{detB}
If $1 \leq l \leq n$ and $2\leq r <n$, then $\det(B_{l,r}) > 0$. 
\end{proposition}

\begin{proof}
The proof is by induction on $l$.
The case for $l=1$ is clear.
For $l=2$, 
$$
\det(B_{2,r})=\det \left[ \begin{array}{rr} 1 & c_{1} \\  -1 & c_{2} \end{array} \right] = c_{2}+c_{1} >0.
$$

Assume $l \geq 3$ and proceed by induction.
By the cofactor expansion along the first row, the determinant of $B_{l,r}$ is
\begin{equation}
\label{cofact}
\det\left[ \begin{array}{rrrccccrc} 
1       & 0      & 0      & 0      & \cdots & \cdots & \cdots & 0       & c_{2} \\
-1      & 1      & 0      &  0      &        &        &        & \vdots  & c_{3} \\
0       & -1     &  1     & 0      &        &        &        &         & \vdots \\ 
\vdots  & 0      & -1     & 1      & \ddots &        &        &         &\vdots  \\
0       &  \vdots   &\ddots  & \ddots & \ddots &        &        &         & \vdots       \\
t_{h}   & 0      &        &        &         &        & \ddots & \vdots  & \vdots \\
0       & t_{h}  & \ddots &        &        &      & \ddots & 0       & c_{l-2} \\
\vdots  & \ddots & \ddots & \ddots &        & \ddots & \ddots &  1      & c_{l-1} \\
0       & \cdots & 0      & t_{h}  & 0      & \cdots & 0      & -1      & c_{l}  
\end{array} \right] + (-1)^{l+1}c_{1}\det(A_{l-1,r}).
\end{equation}

By the induction hypothesis, the first term in~(\ref{cofact}) is positive.
If $r > l-1$, Proposition~\ref{detA} implies that the second term in~(\ref{cofact}) is $(-1)^{l+1}c_{1}(-1)^{l-1}=(-1)^{2l}c_{1} > 0$.
If $2 \leq r \leq l-1$, Proposition~\ref{detA} implies that the second term in~(\ref{cofact}) is
$(-1)^{l+1}c_{1}(-1)^{l-1}a^{0}_{l-1}=c_{1}a^{0}_{l-1} >0.$
Hence, $\det(B_{l,r}) > 0$.
\end{proof}

\medskip

\begin{lemma}
\label{Jacob}
Let $A \in Q(\mathcal{K}_{n,r})$ be of the form~$(\ref{matrix})$. 
When $r<n$, 
the Jacobian $\displaystyle\frac{\partial (\alpha_1, \ldots, \alpha_n)}{\partial (a_1 , \ldots, a_{n-1},b_r)}$
is positive 
at $(1,a^{0}_{2}, \ldots, a^{0}_{n-1},t_h)$.
\end{lemma}

\begin{proof}
By~(\ref{char}), the Jacobian matrix at $(1,a^{0}_{2}, \ldots, a^{0}_{n-1},t_h)$ 
is the block lower triangular matrix
$$J=
\left[ \begin{array}{rrrrcrrcccrrc}
 1    &   0   &  \cdots   &   0   &  \vline     &       &       &        &        &       &       &       &         \\
-1    &   1   &  \ddots   &\vdots &  \vline     &       &       &        &        &    O   &       &       &         \\
      &\ddots &  \ddots   &   0   &  \vline     &       &       &        &       &       &       &       &         \\
0     &       &   -1      &   1   &  \vline     &       &       &        &        &       &       &       &         \\
\noalign{\hrule height 0.3pt}
0     &       &           &   -1  &  \vline     &   1   &   0   & \cdots & \cdots &       & \cdots&  0    & 1        \\
t_{h} &\ddots &           &       &  \vline     &   -1  &   1   &   0    &        &       &           &  0    & 1     \\
      &  t_{h}&  \ddots         &       &  \vline     &   0   &  -1   &   1    &   0    &       &       &  0    & a^{0}_{2}     \\
      &       & \ddots    &  0    &  \vline     &\vdots & \ddots      &\ddots  & \ddots & \ddots&       & \vdots& \vdots     \\
      &       &           & t_{h} &  \vline     &   0   &       &        &\ddots  & \ddots&\ddots & \vdots& \vdots     \\
      &       &           &       &  \vline     & t_{h} & \ddots&        &        & \ddots& \ddots&  0    & \vdots     \\
      &       &           &       &  \vline     &\vdots &\ddots & \ddots &        &\ddots  & -1    &  1    & a^{0}_{n-r-1}     \\
      &       &           &       &  \vline     &  0    &\cdots & t_{h}  &   0    & \cdots&  0     & -1    & a^{0}_{n-r}     
\end{array}\right].
$$
Note that the $(1,1)$-block of $J$ 
has order $r-1$  
and its determinant is $1$.
Since the $t_{h}$ in the last row of $J$ 
is in column $n-r$, 
the $(2,2)$-block of $J$ 
is of the form $B_{n-r+1,r}$.
By Proposition~\ref{detB}, $\det(B_{n-r+1,r}) > 0$.
Thus, the Jacobian at $(1,a^{0}_{2},\ldots, a^{0}_{n-1},t_h)$ is positive.
\end{proof}

\medskip

\begin{theorem}
\label{super}
For $n>r \geq 2$, every superpattern of $\mathcal{K}_{n,r}$ is a SAP.
\end{theorem}

\begin{proof}
Let $A \in Q(\mathcal{K}_{n,r})$ be of the form~(\ref{matrix}).
By Lemma~\ref{NP}, when $a_{1}=\cdots=a_{r-1}=1$, $a_{j}=a^{0}_{j}$ for $j=r, r+1, \ldots, n-1$, and $b_{r}=t_{h}$, 
the resulting matrix is nilpotent.
Moreover, by Lemma~\ref{Jacob}, the Jacobian $\displaystyle\frac{\partial (\alpha_1, \ldots, \alpha_n)}{\partial (a_1 , \ldots, a_{n-1},b_{r})}$ is nonzero 
at $(1,a^{0}_{2}, \ldots, a^{0}_{n-1},t_h)$. 
Thus, Theorem~\ref{theoremNJ} implies that every superpattern of $\mathcal{K}_{n,r}$ is a SAP for $n \geq 3$.
\end{proof}

\medskip

\begin{theorem}
\label{MSAP}
If $n \geq r\geq 2$, then $\mathcal{K}_{n,r}$ is a MSAP.
\end{theorem}

\begin{proof}
As already noted, $\mathcal{K}_{n,n}$ is a MSAP. So assume $n>r$.

We first note that any irreducible subpattern of $\mathcal{K}_{n,r}$ with less than
$2n-1$ nonzero entries is not a SAP~\cite[Theorem 6.2]{BMOD}.

For $i \in \{1, 2, \ldots, n-1\}$, let $\mathcal{S}_{i}$ be the sign pattern obtained by replacing the $(i,i+1)$-entry in $\mathcal{K}_{n,r}$ by 0.
Then $\mathcal{S}_{i}$ is a 2 by 2 block lower triangular matrix.
In particular, the sign of the trace of each of the diagonal blocks of 
$\mathcal{S}_{i}$ is fixed and so each block requires a nonzero 
eigenvalue. Consequently, no subpattern of $\mathcal{S}_{i}$ is 
spectrally arbitrary.

Likewise, if the $(n, n-r+1)$-entry, $(1,1)$-entry, 
or $(n,n)$-entry of $\mathcal{K}_{n,r}$ is replaced by 0, 
the resultant pattern will require a signed eigenvalue.

Next, let $A \in Q(\mathcal{K}_{n,r})$ be of the form~(\ref{matrix}).
Then, by~(\ref{char}), 
$\alpha_{1}=a_{1}-1$, $\alpha_{j}=a_{j}-a_{j-1}$ for $j=2, \ldots, r-1$, $\alpha_{r}=a_{r}+b_{r}-a_{r-1}$, 
$\alpha_{j}=a_{j}+a_{j-r}b_{r}-a_{j-1}$ for $j=r+1, \ldots, n-1$, and $\alpha_{n}=a_{n-r}b_{r}-a_{n-1}$.
Suppose that for some $j \in \{2, \ldots, n-1 \}$, $a_{j}$ is replaced by $0$. 
Then $\alpha_{j+1}$ of the resulting matrix is positive. 
Thus no subpattern of $\mathcal{K}_{n,r}$ is spectrally arbitrary.
\end{proof}

\end{document}